\title{Beauville surfaces and finite groups}

\author{
Yolanda Fuertes\footnote{Supported by grants MTM2006-01859,
MTM2006-14688 and CCG08-UAM/ESP-4145.}\\
Departamento de Matem\'aticas\\
Universidad Aut\'onoma de Madrid\\
Cantoblanco Universidad\\
28049 Madrid\\
Spain\\
{\tt yolanda.fuertes@uam.es}
\and
Gareth A. Jones\\
School of Mathematics\\
University of Southampton\\
Southampton SO17  1BJ\\
U.K.\\
{\tt G.A.Jones@maths.soton.ac.uk}
}

\documentclass[12pt]{article}
\usepackage[latin1]{inputenc}
\usepackage{a4wide}
\usepackage{latexsym}
\usepackage{amsmath}
\usepackage{amssymb}
\usepackage{amsfonts}
\newtheorem{thm}{Theorem}[section]
\newtheorem{lemma}[thm]{Lemma}

\newtheorem{prop}[thm]{Proposition}
\date{}
\begin{document} 

\maketitle

\begin{abstract}
\noindent Extending results of Bauer, Catanese and Grunewald, and of Fuertes and Gonz\'alez-Diez, we show that Beauville surfaces of unmixed type can be obtained from the groups $L_2(q)$ and $SL_2(q)$ for all prime powers $q>5$, and the Suzuki groups $Sz(2^e)$ and the Ree groups $R(3^e)$ for all odd $e\geq 3$. We also show that $L_2(q)$ and $SL_2(q)$ admit strongly real Beauville structures, yielding real Beauville surfaces, for all $q>5$.
\end{abstract}



\noindent{\bf Running head:} Beauville surfaces\\

\section{Introduction}

Algebraic geometers such as Bauer, Catanese and Grunewald~\cite{BCG05, BCG, Cat} have recently initiated the study of Beauville surfaces. These are $2$-dimensional complex algebraic varieties which are rigid, in the sense of admitting no deformations. They are defined over the field $\overline{\bf Q}$ of algebraic numbers, and provide a geometric action of the absolute Galois group ${\rm Gal}\,\overline{\bf Q}/{\bf Q}$. By generalising Beauville's original example~\cite[p.~159]{Bea}, they can be constructed from finite groups acting on suitable pairs of algebraic curves, and here we give some new examples of families of groups which can be used for this purpose.

A {\sl Beauville surface of unmixed type\/} is a compact complex surface $\cal S$ such that
\begin{description}
\item[\rm(a)] $\cal S$ is isogenous to a higher product, that is, ${\cal S}\cong ({\cal C}_1\times {\cal C}_2)/G$ where ${\cal C}_1$ and ${\cal C}_2$ are algebraic curves of genus at least $2$ and $G$ is a finite group acting freely on ${\cal C}_1\times {\cal C}_2$ by holomorphic transformations;
\item[\rm(b)] $G$ acts effectively on each ${\cal C}_i$ so that ${\cal C}_i/G$ is isomorphic to the projective line ${\bf P}^1({\bf C})$ and the covering ${\cal C}_i\to {\cal C}_i/G$ is ramified over at most three points.
\end{description}
(We will not consider the more general situation of a Beauville surface of mixed type, where $G$ contains elements which transpose the two curves ${\cal C}_i$.) Condition~(b) is equivalent to each curve ${\cal C}_i$ admitting a regular {\it dessin\/} in the sense of Grothendieck's theory of {\it dessins d'enfants\/}~\cite{CIW, Gro, Wol}, or equivalently an orientably regular hypermap~\cite{JS}, with $G$ acting as the orientation-preserving automorphism group.

A group $G$ arises in this way if and only if it has generating triples $(x_i, y_i, z_i)$ for $i=1, 2$, of orders $(l_i, m_i, n_i)$, such that
\begin{description}
\item[\rm(1)] $x_iy_iz_i=1$ for each $i=1, 2$,
\item[\rm(2)] $l_i^{-1}+m_i^{-1}+n_i^{-1}<1$ for each $i=1, 2$, and
\item[\rm(3)] no non-identity power of $x_1, y_1$ or $z_1$ is conjugate in $G$ to a power of $x_2, y_2$ or $z_2$.
\end{description}

\noindent We will call such a pair of triples $(x_i, y_i, z_i)$ an {\sl unmixed Beauville structure\/} for $G$, or simply a {\sl Beauville structure\/}. Property~(1) is equivalent to condition~(a), with $x_i$, $y_i$ and $z_i$ representing the ramification over the three points, property~(2) is equivalent to each ${\cal C}_i$ having genus at least $2$ (arising as a smooth quotient of the hyperbolic plane), and property~(3) (which is always satisfied if $l_1m_1n_1$ is coprime to $l_2m_2n_2$) is equivalent to $G$ acting freely on the product.

Bauer, Catanese and Grunewald~\cite{BCG} have made the following conjecture:

\medskip

\centerline{\sl Every non-abelian finite simple group except $A_5$ admits an unmixed Beauville structure.}

\medskip

\noindent They verified that the alternating groups $A_n$ satisfy the conjecture for all sufficiently large $n$, and Fuertes and Gonz\'alez-Diez~\cite{FG} have shown that they do so for all $n\geq 6$. Here we will show that other families of simple groups have this property, namely the groups $L_2(q)$ for prime powers $q>5$, the Suzuki groups $Sz(2^e)$ (extending results for prime $q$ and $e$ in~\cite{BCG}), and the Ree groups $R(3^e)$. We will also show that a family of quasisimple groups (perfect central extensions of simple groups) admit unmixed Beauville structures, namely the groups $SL_2(q)$ for $q>5$, again extending a result for prime $q$ in~\cite{BCG}.  In the case of the groups $L_2(q)$ and $SL_2(q)$, we will show that the Beauville structure can be chosen so that the corresponding Beauville surface is real. We refer to~\cite{BCG, Cat} for background on Beauville surfaces, and to the {\sl ATLAS\/}~\cite{CCNPW} for notation and general information concerning various classes of finite simple and quasisimple groups.

We are very grateful to Prof.~Fabrizio Catanese for pointing out to us that there was some overlap with results recently obtained by Shelly Garion and Matteo Penegini~\cite{GP}. For instance, they use probabilistic methods to show that the Suzuki and Ree groups admit unmixed Beauville structures provided the underlying field is sufficiently large, and they have similar results for other families of simple groups, including $L_3(q)$ and the unitary groups $U_3(q)$, which we have not considered. In the case of the groups $L_2(q)$ they obtain our Theorem~2.2 by using results of Macbeath~\cite{Mac} on generating triples for these groups.

In addition, we thank our colleagues Gabino Gonz\'alez-Diez for valuable advice, and Ernesto Girondo and David Torres for help with computing.

\section{Projective special linear groups $L_2(q)$}

Bauer, Catanese and Grunewald~\cite{BCG} have shown that the simple group $L_2(p)=PSL_2(p)$ admits a Beauville structure for each prime $p>5$ (this fails for $L_2(5)\cong A_5$). We can extend this result to prime powers, but first we need some basic facts about the groups $L_2(q)$; see~\cite[Ch.~XII]{Dic} or~\cite[\S II.8]{Hup} for background.

Let $q=p^e$ for a prime $p$, and let $k=(2,q-1)$. A non-identity element of $L_2(q)$ has order dividing $(q-1)/k$, equal to $p$, or dividing $(q+1)/k$, as it fixes two, one or no points in the projective line ${\bf P}^1(F_q)$ over the field $F_q$. Equivalently, if $t$ is its trace (defined only up to multiplication by $-1$), then $t^2-4$ is respectively a non-zero square, equal to $0$, or a non-square in $F_q$. The group $PGL_2(q)$ contains $L_2(q)$ with index $k$, and its elements have orders dividing $q\pm 1$ or equal to $p$.

Dickson classified the subgroups of $L_2(q)$ (\cite[Ch.~XII]{Dic}, see also~\cite[\S II.8]{Hup}), and from this one can describe the maximal subgroups:

\begin{prop} Any maximal subgroup of $L_2(q)$ has one of the following forms, where $k=(2,q-1)$:

\begin{enumerate}

\item the stabiliser of a point on the projective line ${\bf P}^1(F_q)$, isomorphic to the unique subgroup of order $q(q-1)/k$ in $AGL_1(q)$;

\item a dihedral group of order $2(q\pm 1)/k$;

\item a group isomorphic to $L_2(r)$ where $F_r$ is a maximal subfield of $F_q$;

\item a group isomorphic to $PGL_2(r)$ where $q=r^2$ is a perfect square;

\item a group isomorphic to $A_4$, $S_4$ or $A_5$. \hfill$\square$

\end{enumerate}

\end{prop}

\noindent(Subgroups of types (1) to (4) always exist, but those of type (5) exist only for certain values of $q$, and when they exist they are not always maximal.)

\begin{thm}
For each prime power $q>5$ the group $L_2(q)$ admits a Beauville structure.
\end{thm}

\noindent{\sl Proof.} Let $G=L_2(q)=SL_2(q)/\{\pm I\}$. We will prove this result by choosing elements $X_i, Y_i\in SL_2(q)$ for $i=1, 2$ so that their images $x_i, y_i\in G$ generate $G$, and defining $z_i$ to be the image of $Z_i:=(X_iY_i)^{-1}$, so that $x_iy_iz_i=1$. The orders $l_i, m_i$ and $n_i$ can be controlled by choosing $X_i$, $Y_i$ and $X_iY_i$ to have appropriate traces. Small values of $q$ can be dealt with individually, so we will assume for the moment that $q\geq 11$ if $q$ is odd.

Let
\begin{equation}
X_1=\Big(\,\begin{matrix}0&1\cr -1&a\end{matrix}\,\Big)
\quad{\rm and}\quad
Y_1=\Big(\,\begin{matrix}b&-1\cr 1&0\end{matrix}\,\Big),
\quad{\rm so}\quad
Z_1=\Big(\,\begin{matrix}1&0\cr b-a&1\end{matrix}\,\Big).
\end{equation}
If $q=p^e$ is odd we can choose $a\in F_q$ so that $\pm a$ is the trace of an element of order $(q+1)/2$ in $G$, and then put $b=-a$, so that the elements $x_1$ and $y_1$ of $G$ corresponding to $X_1$ and $Y_1$ have orders $l_1=m_1=(q+1)/2$, while the element $z_1$ corresponding to $Z_1$ has order $n_1=p$.

 By inspecting the maximal subgroups of $G$ in Prop.~2.1 we see that since $(q+1)/2>5$, so that groups of type (5) are excluded,  there is no maximal subgroup containing elements of orders $(q+1)/2$ and $p$. Thus the triple $(x_1, y_1, z_1)$ generates $G$.

If $q=2^e\geq 8$ we can choose distinct values of $a$ and $b$ so that $x_1=X_1$  and $y_1=Y_1$ belong to two distinct conjugacy classes of elements of order $q+1$ in $G=SL_2(q)$ (there are $\phi(q+1)/2\geq 2$ such classes), so that $z_1=Z_1$ has order $n_1=2$; the triple $(x_1, y_1, z_1)$ cannot be contained in a dihedral group, since $x_1$ and $y_1$ have odd order whereas $z$ has order $2$, so again it follows from Prop.~2.1 that it must generate $G$.

We will choose $X_2$ and $Y_2$ in $SL_2(q)$ so that $l_2=m_2=n_2=(q-1)/2$ or $q-1$ as $q$ is odd or even, and hence $l_1m_1n_1$ is coprime to $l_2m_2n_2$. If $q=8$ or $q\geq 13$ it follows from Prop.~2.1 that $x_2$ and $y_2$ generate $G$ provided they have no common fixed point in ${\bf P}^1(F_q)$. Let
\begin{equation}
X_2=\Big(\,\begin{matrix}c&0\cr 0&c^{-1}\end{matrix}\,\Big)
\quad{\rm and}\quad
Y_2=\Big(\,\begin{matrix}x&y\cr z&w\end{matrix}\,\Big),
\quad{\rm so}\quad
Z_2=\Big(\,\begin{matrix}c^{-1}w&-cy\cr -c^{-1}z&cx\end{matrix}\,\Big),
\end{equation}
where $xw-yz=1$. We can choose $c$ so that $l_2$ is as claimed, for instance by taking $c$ to be a primitive root for $F_q$ (i.e.~a generator of the multiplicative group $F_q^*$), or to have order $(q-1)/2$ if this is odd. If we choose $x$ and $w$ so that $x+w=c+c^{-1}$, then ${\rm tr}\,Y_2={\rm tr}\, X_2$ and so $m_2=l_2$. Now
$${\rm tr}\, Z_2=(c-c^{-1})x+(c+c^{-1})c^{-1}\,,$$
with $c-c^{-1}\neq 0$ since $c\neq \pm 1$, so for a fixed $X_2$ there is a bijection between choices of $x$ in $F_q$ and values of ${\rm tr}\,Z_2$. The fixed points of $x_2$ are $0$ and $\infty$. Now $y_2$ fixes these as $y=0$ or $z=0$ respectively, so we need to choose $Y_2$ so that $yz\neq 0$, or equivalently $xw\neq 1$. Since $x+w={\rm tr}\,X_2$ we have $xw=1$ if and only if $\{x, w\}=\{c, c^{-1}\}$, so by letting $x$ avoid these two values we can obtain any value for ${\rm tr}\,Z_2$ except $c^2+c^{-2}$ and $2$. In particular, we can choose $x$ so that ${\rm tr}\,Z_2={\rm tr}\,X_2$, so $n_2=l_2$ as required. If $q=11$ then a triple of elements of order $(q-1)/2=5$ could generate a subgroup $H\cong A_5$; however, a simple calculation within $A_5$ shows that to do so they would need to be conjugate in $H$ and hence in $G$, so a triple such as
\begin{equation}
x_2=\pm\Big(\,\begin{matrix}2&0\cr 0&6\end{matrix}\,\Big),
\quad
y_2=\pm\Big(\,\begin{matrix}0&1\cr -1&-3\end{matrix}\,\Big),
\quad
z_2=\pm\Big(\,\begin{matrix}4&-2\cr -5&0\end{matrix}\,\Big),
\end{equation}
with different traces $\pm 3, \pm 3$ and $\pm 4$, must generate $G$.

This deals with all cases except $q=7$ and $9$. The first is covered by the proof by Bauer, Catanese and Grunewald~\cite{BCG} that $L_2(p)$ admits a Beauville structure for each prime $p>5$. Since $L_2(9)\cong A_6$ the case $q=9$ is covered by the result of Fuertes and Gonz\'alez-Diez~\cite{FG} that the alternating group $A_n$ admits a Beauville structure for each $n\geq 6$.

It is well known and easy to see that the smallest non-abelian finite simple group $L_2(4)\cong L_2(5)\cong A_5$ does not admit a Beauville structure. Each non-identity element of this group has order $2, 3$ or $5$, and any triple consisting of elements of orders $2$ or $3$ would fail to satisfy condition~(2). Any generating triple must therefore contain an element of order $5$, and this violates condition~(3) since all subgroups of order $5$ are conjugate. It is even easier to see that the (non-simple) groups $L_2(3)\cong A_4$ and $L_2(2)\cong S_3$ do not admit Beauville structures. \hfill $\square$

\section{Strongly real Beauville structures.}

A Beauville structure on a group $G$ (unmixed, as before) is {\sl strongly real\/} if there are automorphisms $\alpha_i$ of $G$ for $i=1,2$, differing by an inner automorphism, with each $\alpha_i$ inverting two elements of the triple $(x_i, y_i, z_i)$. This condition implies that the corresponding Beauville surface $\cal S$ is real, that is, there is a biholomorphic map $\sigma: {\cal S}\to\overline{\cal S}$ such that $\sigma^2$ is the identity (see~\cite{BCG} for details). By replacing one triple with its image under that inner automorphism, we may assume that the same automorphism $\alpha$ acts in this way on both triples, and by cyclically permuting the terms of each triple we may assume that it inverts $x_i$ and $y_i$ for $i=1,2$. When $G=L_2(q)$ this implies that $\alpha$ preserves the traces of $x_i, y_i$ and $z_i$: each element of $G$ has the same eigenvalues $\lambda$ and $\lambda^{-1}$ as its inverse, so they have the same trace, and $\alpha$ sends $z_i=y_i^{-1}x_i^{-1}$ to $y_ix_i$ which has the same trace as $x_iy_i$ and hence as $(x_iy_i)^{-1}=z_i$.

Our aim in this section is to show that, with a few small exceptions, each group $L_2(q)$ admits a strongly real Beauville structure. We will do this by adapting the proof of Theorem~2.2. As before, we will assume that $q\geq 8$ if $q$ is even, and $q\geq 11$ if $q$ is odd.

Let $\alpha$ be the automorphism of $G=L_2(q)$ induced by conjugation by the matrix
\begin{equation}
A=\Big(\,\begin{matrix}0&1\cr 1&0\end{matrix}\,\Big)\in GL_2(q).
\end{equation}
An element
\[\Big(\,\begin{matrix}a&b\cr c&d\end{matrix}\,\Big)\in SL_2(q)\]
is inverted by $A$ if and only if $b+c=0$, so $\alpha$ inverts the elements $x_1, y_1$ and $x_2$ used in the proof of Theorem~2.2, and it inverts $y_2$ if $y+z=0$. For instance, the triples used in (1) and (3) for $L_2(11)$ satisfy this condition, so we may assume that $q\geq 13$ if $q$ is odd.

We need to choose $x$ and $w$ as before, but with the additional requirement that $z=-y$, so that $1-xw=-yz=y^2$ must be a square. Let $c$ be a primitive root for $F_q$, and as before let
\begin{equation}
x+w=c+c^{-1},
\end{equation}
so that ${\rm tr}\,Y_2={\rm tr}\,X_2$. Thus $Y_2$, like $X_2$, has eigenvalues $c$ and $c^{-1}$, so these two matrices have order $q-1$. Similarly, if we also let
\begin{equation}
cx+c^{-1}w=c+c^{-1},
\end{equation}
then ${\rm tr}\,Z_2={\rm tr}\,X_2$, so $Z_2$ has order $q-1$. The images $x_2, y_2$ and $z_2$ of $X_2, Y_2$ and $Z_2$ in $G$ therefore have orders $l_2=m_2=n_2=(q-1)/2$ or $q-1$ as $q$ is odd or even. Solving (5) and (6) we find that
\[x=\frac{c^2-c+1-c^{-1}}{c^2-1}=\frac{c+c^{-1}}{c+1}\quad{\rm and}\quad w=\frac{c^2-c+1-c^{-1}}{c-c^{-1}}=cx,\]
so that
\[1-xw=\frac{(c+1)^2-c(c+c^{-1})^2}{(c+1)^2}=\frac{-(c-1)^2(c^2+c+1)}{c(c+1)^2}.\]
In the proof of Theorem~2.2 we required that $xw\neq 1$, so that $x_2$ and $y_2$ have no common fixed points in ${\bf P}^1(F_q)$; here we therefore need $c^3\neq 1$, and this is valid since $q\neq 2, 4$.
We need $1-xw$ to be a square, or equivalently we need
\[-c(c^2+c+1)\]
to be a square. This is always true if $q$ is even, so in this case we can choose $y$ (and $z=-y=y$) so that $x_2$ and $y_2$ are inverted by $\alpha$, as required. The proof of Theorem~2.2 shows that this triple generates $G$.

\medskip

\noindent{\bf Example 3A.} Let $q=8$. We can define $F_8=F_2[t]/(t^3+t+1)$, with $c=t$ generating $F_8^*$. Then $x=t^2$ and $w=t+1$, so $1-xw=t^2+t=(t^2)^2$, and we can take $y=z=t^2$. The matrices
\[X_2=\Big(\,\begin{matrix}t&0\cr 0&t^2+1\end{matrix}\,\Big),
\quad
Y_2=\Big(\,\begin{matrix}t^2&t^2\cr t^2&t+1\end{matrix}\,\Big)
\quad{\rm and}\quad
Z_2=\Big(\,\begin{matrix}t^2&t+1\cr t&t+1\end{matrix}\,\Big)\]
of order $7$ give the required triple in $G=L_2(8)$, with the first two inverted by $\alpha$.

\medskip

We may therefore assume from now on that $q$ is odd, so $q\geq 13$. Since the generator $c$ of $F_q^*$ is now a non-square, we need the element
\[s=-(c^2+c+1)\]
to be a non-square. If this is so, then we can again choose $y$ (and $z=-y$) so that $x_2$ and $y_2$ are inverted by $\alpha$, as required.

\medskip

\noindent{\bf Example 3B.} Let $q=13$. We can choose $c=2$ as a generator for $F_{13}^*$, giving $x=3$ and $w=6$. Then $s=6$ is a non-square, and $1-xw=-4=3^2$, so we can take $y=3$ and $z=-3$. This gives a triple
\[X_2=\Big(\,\begin{matrix}2&0\cr 0&7\end{matrix}\,\Big),
\quad
Y_2=\Big(\,\begin{matrix}3&3\cr -3&6\end{matrix}\,\Big)
\quad{\rm and}\quad
Z_2=\Big(\,\begin{matrix}3&-6\cr -5&6\end{matrix}\,\Big)\]
of matrices of order $12$ in $SL_2(13)$; their images $x_2, y_2$ and $z_2$ in $G=L_2(13)$ have order $6$ and generate $G$, with $x_2$ and $y_2$ inverted by $\alpha$.

\medskip

We may therefore assume that $q>13$ and $s$ is a square. Instead of (6), let us impose the condition that
\begin{equation}
cx+c^{-1}w=-c-c^{-1},
\end{equation}
so that $Z_2$ has eigenvalues $-c$ and $-c^{-1}$, and hence has order $q-1$ or $(q-1)/2$ as $q\equiv 1$ or $-1$ mod~$(4)$. Thus $x_2, y_2$ and $z_2$ have orders $l_2=m_2=n_2=(q-1)/2$ as before. On solving (5) and (7) we find that
\[x=\frac{-c^2-c-1-c^{-1}}{c^2-1}=\frac{c+c^{-1}}{1-c}\quad{\rm and}\quad w=\frac{c^2+c+1+c^{-1}}{c-c^{-1}}=-cx,\]
so that
\[1-xw=\frac{(1-c)^2+c(c+c^{-1})^2}{(1-c)^2}=\frac{(c+1)^2(c^2-c+1)}{c(1-c)^2}.\]
The condition $xw\neq 1$ is satisfied provided $c^3\neq -1$, and this is valid since $q>7$.
In this case, in order for $1-xw$ to be a square we need the element
\[t=c^2-c+1\]
to be a non-square. If $t$ is a square then the element
\[st=-(c^4+c^2+1),\]
as a product of two squares, is also a square. In this case we can go back and replace $c$ with $c^2$ in our original choice of $x$ and $w$ (equations~(2), (5) and (6)), so that $X_2$, $Y_2$ and $Z_2$ have eigenvalues $c^2$ and $c^{-2}$ and hence have order $(q-1)/2$. This gives a triple $(x_2, y_2, z_2)$ in $G$ which have orders $l_2=m_2=n_2=(q-1)/2$ or $(q-1)/4$ as $q\equiv -1$ or $1$ mod~$(4)$ respectively. In order that $xw\neq 1$ we now require $c^6\neq 1$, valid since $q>7$. In order for $1-xw$ to be a square we require
\[-c^2(c^4+c^2+1)\]
to be a square, and this is true since $st$ is a square. We can therefore choose $y$ and $z=-y$ as before, giving the required triple.

\medskip

\noindent{\bf Example 3C.} This last situation can arise. For instance, if $q=37$ and we choose $c=2$ as a generator for $F_{37}^*$, then $s=-7$ and $t=3$ are both squares, namely of $17$ and $15$. If we replace $c=2$ with $2^2=4$ then equations~(5) and (6) give $x=-1$ and $w=-4$, so $1-xy=-3=16^2$; taking $y=16$ and $z=-16$ gives a triple
\[X_2=\Big(\,\begin{matrix}4&0\cr 0&-9\end{matrix}\,\Big),
\quad
Y_2=\Big(\,\begin{matrix}-1&16\cr -16&-4\end{matrix}\,\Big)
\quad{\rm and}\quad
Z_2=\Big(\,\begin{matrix}-1&10\cr 4&-4\end{matrix}\,\Big)\]
of elements of order $18$ in $SL_2(37)$, and hence a triple $(x_2, y_2, z_2)$ of elements of order $9$ in $G=L_2(37)$, with $x_2$ and $y_2$ inverted by $\alpha$.

\medskip

It remains for us to show that in this last situation, if $q\equiv 1$ mod~$(4)$ then the triple $(x_2, y_2, z_2)$ of elements of order $(q-1)/4$ generates $G$. We have ensured that $xw\neq 1$, so $x_2$ and $y_2$ have no common fixed points in ${\bf P}^1(F_q)$ and cannot therefore be contained in maximal subgroups of type~(1) or (2) in Prop.~2.1. If $F_r$ is a proper  subfield of $F_q$ and $L_2(r)$ contains elements of order $(q-1)/4$, then $(r^2-1)/4\leq (q-1)/4\leq (r+1)/2$, giving $r\leq 3$ and hence $q\leq 9$, against our assumption. Similarly, if $q=r^2$ and $PGL_2(r)$ contains elements of order $(q-1)/4$ then $(r^2-1)/4\leq r+1$, so $r\leq 5$ and hence $q\leq 25$; since $q$ is an odd square, greater than $13$, we must have $q=25$ with $PGL_2(r)\cong S_5$, whereas elements of order $6$ in $S_5$ are all odd and hence cannot satisfy $x_2y_2z_2=1$. This leaves only subgroups isomorphic to $A_4$, $S_4$ or $A_5$ as possible maximal subgroup containing the triple. If $q\geq 25$ then since these have no elements of order $6$ or higher, we are done. The only remaining possibility is that $q=17$ and $x_2, y_2$ and $z_2$ correspond to elements of order $4$ in $S_4$, again impossible since they would all be odd.

This deals with $q=8$ and all $q\geq 11$. In the case $q=7$ the triple of matrices
\[X_1 = \Bigl(\begin{matrix} 0 & 1 \\ -1 & 3 \end{matrix}\Bigr),
\quad Y_1 = \Bigl(\begin{matrix} -2 & 2 \\ -2 & -2 \end{matrix}\Bigr)
\quad{\rm and}\quad Z_1 = \Bigl(\begin{matrix} -2 & -2 \\ 3 & -1 \end{matrix}\Bigr)\]
of order $8$ in $SL_2(7)$ have images $x_1, y_1$ and $z_1$ of order $4$ in $G=L_2(7)$, while the triple\[X_2 = \Bigl(\begin{matrix} 0 & 1 \\ -1 & 2 \end{matrix}\Bigr),
\quad Y_2 = \Bigl(\begin{matrix} 0 & -1 \\ 1 & 2 \end{matrix}\Bigr)
\quad{\rm and}\quad Z_2 = \Bigl(\begin{matrix} -2 & -2 \\ -2 & 1 \end{matrix}\Bigr)\]
of order $7, 7$ and $3$ have images $x_2, y_2$ and $z_2$ of these same orders in $G$. Each triple $(x_i, y_i, z_i)$ generates $G$: as in the preceding case $q=17$ the first triple cannot generate a subgroup isomorphic to $S_4$, and the second cannot lie in a point stabiliser since $x_2$ and $y_2$ fix different points in ${\bf P}^1(F_7)$. Since $x_i$ and $y_i$ are inverted by $\alpha$ for each $i$, the resulting Beauville structure on $G$ is strongly real.

For $q=9$ we can take $F_9=F_2[t]/(t^2+1)$. The triple of matrices
\[X_1 = \Bigl(\begin{matrix} t+1 & 0 \\ 0 & t-1 \end{matrix}\Bigr),
\quad Y_1 = \Bigl(\begin{matrix} -t+1 & t \\ -t+1 & -1 \end{matrix}\Bigr)
\quad{\rm and}\quad Z_1 = \Bigl(\begin{matrix} -t+1 & -t+1 \\ t & -1 \end{matrix}\Bigr)\]
of order $8$ in $SL_2(9)$ have images $x_1, y_1$ and $z_1$ of order $4$ in $G=L_2(9)$. As in earlier cases, since  $x_1, y_1$ and $z_1$ fix different points in ${\bf P}^1(F_9)$, and cannot generate a subgroup isomorphic to $S_4$, they generate $G$. The triple
\[X_2 = \Bigl(\begin{matrix} 1 & t+1 \\ t & t \end{matrix}\Bigr),
\quad Y_2 = \Bigl(\begin{matrix} t & t+1 \\ t & 1 \end{matrix}\Bigr)
\quad{\rm and}\quad Z_2 = \Bigl(\begin{matrix} -t -1 & t+1 \\ -1 & -t-1 \end{matrix}\Bigr)\]
in $SL_2(9)$ all have order $5$, and have images $x_2, y_2$ and $z_2$ of the same order in $G$. By Proposition~2.1, if $x_2, y_2$ and $z_2$ generate a proper subgroup $H<G$ then $H\cong C_5$ or $H\cong A_5\cong L_2(5)$, and since they do not commute we must have $H\cong L_2(5)$; it follows then that $X_2, Y_2$ and $Z_2$ generate a subgroup isomorphic to $SL_2(5)$, the only covering group of $H$ in $SL_2(9)$. However, $SL_2(5)$ is not generated by any triple $(X, Y, Z)$ of type $(5,5,5)$: elements of order $5$ in $SL_2(5)$ have trace $2$, so without loss of generality (by applying an automorphism) we may assume that 
\[X = \Bigl(\begin{matrix} 1 & 1 \\ 0 & 1 \end{matrix}\Bigr)
\quad{\rm and}\quad
\quad Y = \Bigl(\begin{matrix} x & y \\ z & 2-x \end{matrix}\Bigr)\]
for some $x, y, z\in F_5$; then ${\rm tr}\,Z={\rm tr}\,XY=z+2$, so $z=0$ and hence $X$ and $Y$ generate a proper subgroup of $SL_2(5)$. This shows that the triple $(x_2, y_2, z_2)$ generates $G$. For each $i=1, 2$ the matrices $X_i$ and $Y_i$ are inverted by conjugation by the matrix
\[B=\Bigl(\begin{matrix} 0 & 1 \\ t+1 & 0 \end{matrix}\Bigr)\in GL_2(9),\]
so $x_i$ and $y_i$ are inverted by the corresponding automorphism of $G$. Thus the triples $(x_i, y_i, z_i)$ form a strongly real Beauville structure on $G$.

Since $L_2(q)$ admits no Beauville structures for $q\leq 5$, we have therefore proved:

\begin{thm}
The group $L_2(q)$ admits a strongly real Beauville structure if and only if $q>5$. \hfill$\square$
\end{thm}

This result provides partial evidence for a more ambitious conjecture of Bauer, Catanese and Grunewald in~\cite{BCG05} that all but finitely many non-abelian finite simple groups admit a strongly real unmixed Beauville structure.

\section{Lifting Beauville structures}

Bauer, Catanese and Grunewald~\cite{BCG} have shown that the group $SL_2(p)$ admits a Beauville structure for each prime $p>5$. Again we can extend this result to prime powers, but first we need some preparatory results.

When proving that a composite group $G$, such as $SL_2(q)$ for odd $q$, admits a Beauville structure, it is tempting to look for such a structure in the quotient $G/N$ by some normal subgroup $N\neq 1$ of $G$, and to try to lift this back to $G$. However, a triple that generates $G/N$ need not lift back to a triple generating $G$, and even if it does, condition~(1) may not be satisfied. If these difficulties can be overcome, then there is no problem with condition~(2), since lifting cannot decrease the orders of elements. However, condition~(3) may be troublesome, since cyclic subgroups which have trivial intersection in $G/N$ need not lift back to subgroups with this property in $G$. The following example is instructive.

\medskip

\noindent{\bf Example 4A.} Let $G$ be the metacyclic group of order $p^3$ with presentation
\[\langle a, b \mid a^{p^2}=b^p=1,\; a^b=a^{p+1}\rangle,\]
where $p$ is prime. This has a normal subgroup $N=G'=Z(G)=\langle a^p\rangle\cong C_p$ with $G/N\cong C_p\times C_p$. If $p\geq 5$ then $G/N$ admits a Beauville structure. (This is because it has $p+1\geq 6$ subgroups of order $p$, which is enough to allow the choice of two suitable triples; the corresponding curves ${\cal C}_i$ are the Fermat curves ${\cal F}_p$ of genus $(p-1)(p-2)/2$, given in homogeneous coordinates by $x^p+y^p+z^p=0$.) However, if $p\geq 3$ then all elements of $G$ have order $p^2$, apart from those in $\langle a^p, b\rangle$. It follows that any generating triple must contain at least one (in fact two) elements $g$ of order $p^2$. Thus $\langle g\rangle$ contains $\langle g^p\rangle=N$, so no two triples can satisfy condition~(3), and hence $G$ does not admit a Beauville structure.

\begin{lemma} If $G$ is a perfect group, $N$ is a central subgroup of $G$, and $S$ is a subset of $G$ such that the image of $S$ in $G/N$ generates $G/N$, then $S$ generates $G$.
\end{lemma}

\noindent{\sl Proof.} Let $H$ be the subgroup of $G$ generated by $S$. Then $HN=G$, so $H$ is a normal subgroup of $G$ since it is normalised by itself and by the central subgroup $N$. Now $G/H=HN/H\cong N/(N\cap H)$, so $G/H$ is abelian since $N$ is. However, $G$ is perfect, so $H=G$. \hfill$\square$

\medskip

This shows that in a quasisimple group (a perfect central extension of a simple group), any subset which maps onto a generating set for the simple quotient must generate the whole group. In particular, a subset of $SL_2(q)$ generates $SL_2(q)$ if and only if its image in $L_2(q)$ generates $L_2(q)$.

\medskip

If $G$ is a group with a normal subgroup $N$, we say that an element $g$ of $G$ is faithfully represented in $G/N$ if $\langle g\rangle\cap N=1$, or equivalently the order of $g$ in $G$ is the same as that of its image in $G/N$. We say that a triple in $G$ is faithfully represented in $G/N$ if each of its elements is faithfully represented in $G/N$.

\begin{lemma} Let $G$ have generating triples $(x_i, y_i, z_i)$ with $x_iy_iz_i=1$ for $i=1, 2$, and a normal subgroup $N$ such that at least one of these triples is faithfully represented in $G/N$. If the images of these triples correspond to a Beauville structure for $G/N$, then these triples  correspond to a Beauville structure for $G$.
\end{lemma}

\noindent{\sl Proof.} Without loss of generality we may assume that $(x_1, y_1, z_1)$ is faithfully represented in $G/N$. Now suppose, without loss of generality, that $x_1^j$ is conjugate in $G$ to a power of $x_2, y_2$ or $z_2$. Then the image of $x_1^j$ in $G/N$ is conjugate in $G/N$ to a power of the image of $x_2, y_2$ or $z_2$. The Beauville property for $G/N$ implies that this image must be the identity, so $x_1^j\in N$ and hence $x_1^j=1$. \hfill$\square$

\section{Special linear groups $SL_2(q)$}

Here we will apply the results of \S 4 to the groups $SL_2(q)$.

\begin{thm}
The group $SL_2(q)$ admits a Beauville structure if and only if $q>5$. In all such cases it admits a strongly real Beauville structure.
\end{thm}

\noindent{\sl Proof.} Let $G=SL_2(q)$ and $\overline G=L_2(q)$. If $q=2^e$ then $G=\overline G$, so Theorems~2.2 and 3.1 give the result. We may therefore assume that $q$ is odd, so that $G$ is a double covering of $\overline G$. We assume first that $q\geq 13$, using Lemmas~4.1 and 4.2 to deduce the result from the methods of proof of Theorems~2.2 and 3.1; smaller values of $q$ are dealt with later by separate arguments. Specifically, we use Lemma~4.1 to show that generating triples $(x_i, y_i, z_i)$ for $\overline G$ lift back to triples $(X_i, Y_i, Z_i)$ which generate $G$. This allows us to use Lemma~4.2 to lift Beauville structures from $\overline G$ to $G$; in order to satisfy the hypotheses of Lemma~4.2 we choose the matrices $X_i, Y_i$ and $Z_i$ so that in each case one of the two triples $(X_i, Y_i, Z_i)$ consists of elements of odd order, and is therefore faithfully represented in $\overline G$. As in the case of $L_2(q)$, in order to obtain a strongly real Beauville structure we choose $X_i$ and $Y_i$ to be inverted by conjugation by the matrix $A$ in (4) for $i=1, 2$.

\smallskip

\noindent{\bf Case 1.} Suppose first that $q\equiv 1$ mod~$(4)$, with $q\geq 13$. Since $q>5$ there exist elements $u$ and $v\neq u^{\pm 1}$ of order $(q+1)/2$ in $F_{q^2}$, so $a:=u+u^{-1}$ and $b:=v+v^{-1}$ are distinct elements of $F_q$. Using these values of $a$ and $b$ we define $X_1$, $Y_1$ and $Z_1$ as in equation~(1). Since $X_1$ and $Y_1$ have eigenvalues $u^{\pm 1}$ and $v^{\pm 1}$, they have order $(q+1)/2$, while $Z_1$ has order $p$. These are all odd, so the triple $(X_1, Y_1, Z_1)$ is faithfully represented in $\overline G$. Since $(q+1)/2>5$ it follows from Proposition~2.1 that $\overline G$ is generated by the image of this triple, so Lemma~4.1 implies that $(X_1, Y_1, Z_1)$ generates $G$. A similar argument shows that the triple $(X_2, Y_2, Z_2)$ defined in the proof of Theorem~3.1 also generates $G$. By their construction, these three matrices all have orders $q-1$ or $(q-1)/2$, coprime to the orders of $X_1$, $Y_1$ and $Z_1$. These two triples therefore form a Beauville structure for $G$. Moreover, since conjugation by $A$ inverts $X_i$ and $Y_i$ for $i=1, 2$, this structure is strongly real. This argument fails when $q\leq 9$ since the chosen triples need not generate $G$; we will deal with this case later.

\smallskip

\noindent{\bf Case 2.} Now suppose that $q\equiv -1$ mod~$(4)$, with $q>11$. We choose $X_1$, $Y_1$ and $Z_1$ as in case 1, but with $u$ and $v$ now of order $q+1$; thus $X_1$ and $Y_1$ are inverted by $A$ and have order $q+1$, which is even, while $Z_1$ again has order $p$. The images of $X_1$ and $Y_1$ in $\overline G$ have order $(q+1)/2>5$, while that of $Z_1$ has order $p$, so it again follows that the triple $(X_1, Y_1, Z_1)$ generates $G$. We now need a triple $(X_2, Y_2, Z_2)$ consisting of elements of odd order dividing $q-1$. If we ignore the requirement that $X_2$ and $Y_2$ should be inverted by $A$, then it is easy to construct a Beauville structure for $G$: since $q>11$ we can choose $c\in F_q$ as in the proof of Theorem~2.2 so that the matrices $X_2$, $Y_2$ and $Z_2$ in $(2)$ have odd order $(q-1)/2$ and generate $G$ (since their images generate $\overline G$). 

The matrix $X_2$ in $(2)$ is inverted by $A$ for any choice of $c$, but in order to construct a strongly real Beauville structure we also need $Y_2$ to be inverted by $A$, and this happens if and only if $z=-y$, as in the proof of Theorem~3.1. We therefore define $X_2, Y_2$ and $Z_2$ by
\begin{equation}
X_2=\Big(\,\begin{matrix}c&0\cr 0&c^{-1}\end{matrix}\,\Big)
\quad{\rm and}\quad
Y_2=\Big(\,\begin{matrix}x&y\cr -y&w\end{matrix}\,\Big),
\quad{\rm so}\quad
Z_2=\Big(\,\begin{matrix}c^{-1}w&-cy\cr c^{-1}y&cx\end{matrix}\,\Big).
\end{equation}
We again use equations~(5) and (6), so that $x+w=cx+c^{-1}w=c+c^{-1}$, but now taking $c=-d$ where $d$ is a primitive root for $F_q$, so that $c$ has order $(q-1)/2$. As before, we have
\[1-xw=\frac{-(c-1)^2(c^2+c+1)}{c(c+1)^2}.\]
Since $c$ is now a square in $F_q$ whereas $-1$ is not, in order for $1-xw$ to be a square we require the element $c^2+c+1=d^2-d+1$ to be a non-square. If this is the case, we can find a triple $(X_2, Y_2, Z_2)$ of elements which have odd order $(q-1)/2$, since they have eigenvalues $c^{\pm 1}$; they generate $G$, and $X_2$ and $Y_2$ are inverted by $A$, so we have a strongly real Beauville structure.

\smallskip

\noindent{\bf Case 3.} Now suppose that $11<q\equiv -1$ mod~$(4)$ as before, and $d^2-d+1$ is a square for each primitive root $d\in F_q$. (This happens if $p=3$, for instance, since then $d^2-d+1=(d+1)^2$.) In (8) we put $c=-d$ as in case~(2), but now with $x+w=c^2+c^{-2}$ and $cx+c^{-1}w=c^2+c^{-2}$. Then $X_2$ has eigenvalues $-d^{\pm 1}$, while $Y_2$ and $Z_2$ have eigenvalues $c^{\pm 2}=d^{\pm 2}$, so they all have order $(q-1)/2>5$. Solving the two equations for $x$ and $w$ we obtain
\[x=\frac{d^4+1}{d^2(1-d)}
\quad{\rm and}\quad
w=\frac{d^4+1}{d(d-1)},\]
so
\[1-xw=\frac{(d+1)^2(d^2-d+1)(d^4-d^3+d^2-d+1)}{d^3(d-1)^2}.\]
This is a non-zero square, giving us a strongly real Beauville structure on $G$, provided $(d^2-d+1)(d^4-d^3+d^2-d+1)$ is a non-square. Since $d^2-d+1$ is a square, non-zero since $q>7$, this is equivalent to $d^4-d^3+d^2-d+1$ being a non-square.

\smallskip

\noindent{\bf Case 4.} Now suppose that $11<q\equiv -1$ mod~$(4)$ as before, and that $d^2-d+1$ and $d^4-d^3+d^2-d+1$ are both squares for each primitive root $d\in F_q$. In (8) we put $c=-d$ again, but now with $x+w=c^3+c^{-3}$ and $cx+c^{-1}w=c+c^{-1}$. Then $X_2$ and $Z_2$ have order $(q-1)/2$, while $Y_2$, with eigenvalues $c^{\pm 3}=-d^{\pm 3}$, has order $(q-1)/6$ or $(q-1)/2$ as $q\equiv 1$ mod~$(3)$ or not. We have
\[
x=\frac{d^5+d^2-d+1}{d^3(d-1)}
\quad{\rm and}\quad
w=\frac{d^5-d^4+d^3+1}{d(1-d)},
\]
so
\[
1-xw=\frac{(d^4-d^3+d^2-d+1)(d^3+1)^2}{d^4(d-1)^2}.
\]
This is a square, non-zero since $q>11$, so we obtain a strongly real Beauville structure.

\smallskip

Having dealt with all the prime powers $q\geq 13$, we now consider small values of $q$.

\smallskip

\noindent{\bf Case 5.} Let $q=11$. The arguments in cases~(2), (3) and (4) do not apply to $G=SL_2(11)$, and the strongly real Beauville structure for $\overline G=L_2(11)$ given by the triples~(1) and (3) does not lift back to a Beauville structure for $G$ since at least one of the elements $x_2, y_2$ and $z_2$ of order $5$ in (3) must lift back to an element of order $10$, violating condition~(3). Instead, consider the triples
\[
X_1=\left( \begin{matrix}
0 & 1 \\ -1 & 5
\end{matrix}\right),
\quad
Y_1=\left( \begin{matrix}
0 & 1\\ -1 & -5
\end{matrix}\right)
\quad{\rm and}\quad
Z_1=\left( \begin{matrix}
-4 & 5 \\ 5 & -1
\end{matrix}\right),\]
all of order $12$, and
\[
X_2=\left( \begin{matrix}
0 & 1 \\ -1 & -4
\end{matrix}\right),
\quad
Y_2=\left( \begin{matrix}
3 & -1\\ 1 & 0
\end{matrix}\right)
\quad{\rm and}\quad
Z_2=\left( \begin{matrix}
1 & 0 \\ -4 & 1
\end{matrix}\right),\]
of orders $5, 5$ and $11$. The images of $X_1$ and $Y_1$ in $\overline G$ have order $6$ and do not commute, so it follows from Proposition~2.1 and Lemma~4.1 that the triple $(X_1, Y_1, Z_1)$ generates $G$. The image of $X_1$ fixes $3$ and $4$ in ${\bf P}^1(F_{11})$, whereas the image of $Z_2$ fixes only $0$, so the triple $(X_2, Y_2, Z_2)$ also generates $G$. In each case $X_i$ and $Y_i$ are inverted by $A$, so the resulting Beauville structure on $G$ is strongly real.

\smallskip

\noindent{\bf Case 6.} For $G=SL_2(9)$ we can use the triples $(X_i, Y_i, Z_i)$ defined in the case $q=9$ of the proof of Theorem~3.1: they generate $G$ by Lemma~4.1, and each matrix $X_i$ or $Y_i$ is inverted by conjugation by $B$, so they form a strongly real Beauville structure on $G$.

\smallskip

\noindent{\bf Case 7.}  Bauer, Catanese and Grunewald~\cite{BCG} have shown that $SL_2(p)$ admits a Beauville structure for each prime $p>5$, so this applies to $SL_2(7)$. In fact, the matrices $X_i, Y_i$ and $Z_i$ defined in the case $q=7$ of the proof of Theorem~3.1 give a strongly real Beauville structure on this group.

\smallskip

\noindent{\bf Case 8.}  As in the case of $L_2(5)$, there is no Beauville structure on $SL_2(5)$: any generating triple for this group must contain an element of order $5$ or $10$, since it maps onto a generating triple for $L_2(5)$; however $SL_2(5)$ has a single conjugacy class of cyclic subgroups of order $10$, and these contain all the elements of order $5$, so any two generating triples must violate condition~(3). A similar argument, based on elements of order $3$, gives the same result for $SL_2(3)$. \hfill$\square$

\medskip

\noindent{\bf Example 5A.} As an illustration of Theorem~5.1(b) with $11<q\equiv -1$ mod~$(4)$, suppose that $q$ is prime, that $d=2$ is a primitive root for $F_q$, so that $q\equiv \pm 3$ mod`$(8)$, and that $q\equiv 1$ mod~$(3)$, giving $q\equiv 19$ mod~$(24)$. The element $d^2-d+1=3$ is a non-square mod~$(q)$ by quadratic reciprocity, since $q\equiv 3\equiv -1$ mod~$(4)$ and $q$ is a square mod~$(3)$; we can therefore take $c=-d=-2$ in case~(2) of the above proof to obtain a strongly real Beauville structure on $SL_2(q)$. (E.~Artin conjectured that the set of primes for which $2$ is a primitive root has asymptotic density
\[\prod_{p\;{\rm prime}}\Bigl(1-\frac{1}{p(p-1)}\Bigr)=0.3739558136\ldots\]
in the set of all primes; this is still unproved.) For instance, if $q=19$ then $2$ is a primitive root; putting $c=-2$ in case~(2) gives $x=-7$ and $w=-5$, so $1-xw=4$, which is a square; taking $y=2$, so that $xw+y^2=1$, we obtain a triple
\[
X_2=\left( \begin{matrix}
-2 & 0 \\ 0 & 9
\end{matrix}\right),
\quad
Y_2=\left( \begin{matrix}
-7 & 2\\ -2 & 5
\end{matrix}\right),
\quad
Z_2=\left( \begin{matrix}
-7 & 4 \\-1 & -5
\end{matrix}\right),\]
of elements of order $9$, forming part of a strongly real Beauville structure on $SL_2(19)$. The other triple $(X_1, Y_1, Z_1)$, given by (1), consists of elements of orders $20, 20$ and $19$.

\medskip

More generally, if $11<q\equiv -1$ mod~$(4)$ then in order to produce a specific strongly real Beauville structure for $SL_2(q)$ we need to know whether either of $d^2-d+1$ and $d^4-d^3+d^2-d+1$ is a square in $F_q$ for a given primitive root $d$, so that we can apply the construction in case~(2), (3) or (4). Quadratic reciprocity deals with this when $q$ is prime, but if $e>1$ then we need Dedekind's generalisation of this law to all finite fields.

If $q=p^e$ with $p$ prime, then $F_q$ can be represented as $F_p[t]/(f(t))$ where $f(t)$ is an irreducible polynomial of degree $e$ in $F_p[t]$. It is convenient to take $f(t)$ to be a primitive polynomial, that is, the minimal polynomial of a primitive root of $F_q$, so that the primitive roots are the powers $t^i$ with $(i,q-1)=1$. The elements of $F_q$ are uniquely represented as the polynomials $g(t)\in F_p[t]$ of degree less than $e$. In testing whether $g(t)$ is a square in $F_q$, one may assume that $q$ is odd (since every element is a square if $q=2^e$), and that $g(t)$ is monic:  if $e$ is even then every constant $a\in F_p$ is a square in $F_q$, and if $e$ is odd then $a$ is a square in $F_q$ if and only if it is a square in $F_p$, which can be tested by classical quadratic reciprocity.

Dedekind's extension of quadratic reciprocity is as follows~\cite{BS}. Let $f(t)$ be a non-constant irreducible monic polynomial over a field $F$ of odd order. Given any polynomial $g(t)\in F[t]$ we define $(g/f)$ to be $+1$, $-1$ or $0$ as $g(t)$ represents a non-zero square, a non-square, or $0$ in the field $F[t]/(f(t))$. More generally, if $f(t)$ is a product of non-constant irreducible monic polynomials $f_i(t)\in F[t]$ we define $(g/f)=\prod_i(g/f_i)$. Dedekind showed that
\begin{equation}
\Bigl(\frac{g}{f}\Bigr)\Bigl(\frac{f}{g}\Bigr)=(-1)^{\deg(f)\deg(g)(|F|-1)/2}.
\end{equation}
In our case we will use this with $F=F_p$ where $p\equiv -1$ mod~$(4)$, so (9) simplifies to
\begin{equation}
\Bigl(\frac{g}{f}\Bigr)\Bigl(\frac{f}{g}\Bigr)=(-1)^{\deg(f)\deg(g)}.
\end{equation}

\medskip

\noindent{\bf Example 5B.} Let $q=3^3=27$. The polynomial $f(t)=t^3-t+1\in F_3[t]$ is primitive, so $F_{27}=F_3[t]/(f(t))$, and we can take $d=t$ as a primitive root. As in all cases where $p=3$, we have $d^2-d+1=(d+1)^2$, a square. Using $t^3=t-1$ and $t^4=t^2-t$ we find that $d^4-d^3+d^2-d+1=-t^2-1$. Then
\[\Bigl(\frac{-t^2-1}{f(t)}\Bigr)
=-\Bigl(\frac{t^2+1}{f(t)}\Bigr)
=-\Bigl(\frac{f(t)}{t^2+1}\Bigr)
=-\Bigl(\frac{t+1}{t^2+1}\Bigr)
=-\Bigl(\frac{t^2+1}{t+1}\Bigr)
=-\Bigl(\frac{2}{3}\Bigr)
=1,
\]
where we have used $f(t)=t(t^2+1)+t+1$, and $(2/3)$ is the Legendre symbol. Thus $d^4-d^3+d^2-d+1$ is a non-zero square (of $t^2-t$, in fact), so we use the construction in case~(4) of the proof.

Putting $d=t$ in case~(4) we find that $x=0$ and $w=t^2+1$, so $1-xw=1$. We can therefore put $y=1$, giving a triple
\[
X_2=\Bigl(\begin{matrix} -t & 0 \\ 0 & -t^{-1} \end{matrix}\Bigr),
\quad
Y_2=\Bigl(\begin{matrix} 0 & 1 \\ -1 & t^2+1 \end{matrix}\Bigr)
\quad{\rm and}\quad
Z_2=\Bigl(\begin{matrix} t^2-t-1 & t \\ t^2-1 & 0 \end{matrix}\Bigr),
\]
all of order $(q-1)/2=13$. As usual, (1) gives the other triple $(X_1, Y_1, Z_1)$.

\medskip

In the case where $11<q\equiv -1$ mod~$(4)$, since exactly half of the elements of $F_q^*$ are squares, one might expect that on average, $d^2-d+1$ should be a non-square (equivalently $d-1+d^{-1}$ should be a square) for about half of the $\phi(q-1)/2$ inverse pairs $d^{\pm 1}$ of primitive roots in $F_q$. The existence of at least one such pair would allow us to use the construction in case~(2) for a strongly real Beauville structure on $SL_2(q)$. As $q$ becomes large, so does $\phi(q-1)/2$, so it seems increasingly likely that such a pair should exist. As supporting evidence, Table~1 shows the primes $q\equiv -1$ mod~$(4)$ from $11$ to $103$, with a primitive root $d$ (expressed as the least possible power of the smallest primitive root) such that $r:=d-1+d^{-1}$ is a quadratic residue mod~$(q)$, that is, a square in $F_q^*$.

\bigskip

\noindent\begin{tabular}{||c||c|c|c|c|c|c|c|c|c|c|c|c||} \hline
$q$ & 11 & 19 & 23 & 31 & 43 & 47 & 59 & 67 & 71 & 79 & 83 & 103 \\ \hline
$d$ & $2^3=8$ & 2 & 5 & $3^{13}=24$ & 3 & $5^3=31$ & $2^7=10$ & 2 & $7^7=14$ & 3 & $2^3=8$ & 5\\ \hline
$d^{-1}$ & 7 & 10 & 14 & 22 & 29 & 44 & 6 & 34 & 66 & 53 & 52 & 62\\ \hline
$r$ & 3 & 11 & 18 & 14 & 31 & 27 & 15 & 35 & 8 & 55 & 59 & 66 \\  \hline
\end{tabular}

\medskip

\centerline{Table 1}


\bigskip

On the basis of this we conjecture that  $F_q$ possesses such a primitive root for every prime $q\equiv -1$ mod~$(4)$, $q\geq 11$. However, if $q=3^e$ then $r=(d+1)^2/d$ is never a square, and we are forced to use the construction in case~(3) or case~(4).

\medskip

Theorem~5.1 suggests the following variation of the conjecture in \S 1:

\medskip

\centerline{\sl Does every finite quasisimple group except $L_2(5)$ and $SL_2(5)$ admit a Beauville structure?}

\medskip

\noindent Similarly, Theorem~5.1 raises the question of which other quasisimple groups admit strongly real Beauville structures.

\section{Suzuki groups and Ree groups}

We now return to the original conjecture concerning finite simple groups. The Suzuki group $Sz(q)={}^2B_2(q)$ is a simple group of order $q^2(q^2+1)(q-1)$, where $q=2^e$ for some odd $e\geq 3$. Bauer, Catanese and Grunewald~\cite{BCG} have shown that $Sz(2^e)$ admits a Beauville structure whenever $e$ is prime. We can extend this result to all Suzuki groups. First we need a general result which allows us to count triples of a given type in a finite group:

\begin{prop}
{\rm\cite[\S 7.2]{Ser}} If $\cal X$, $\cal Y$ and $\cal Z$ are conjugacy classes in any finite group $G$, then the number $N({\cal X},{\cal Y},{\cal Z})$ of solutions of $xyz=1$ with $x\in{\cal X}$, $y\in{\cal Y}$ and $z\in{\cal Z}$ is given by
$$N({\cal X},{\cal Y},{\cal Z})=\frac{|{\cal X}|.|{\cal Y}|.|{\cal Z}|}{|G|}\sum_{\chi}\frac{\chi(x)\chi(y)\chi(z)}{\chi(1)},$$
where $\chi$ ranges over the irreducible complex characters of $G$. \hfill$\square$
\end{prop}

\begin{thm}
The Suzuki group $Sz(2^e)$ admits a Beauville structure for each odd $e\geq 3$.
\end{thm}

\noindent{\sl Proof.} Suzuki~\cite{Suz} showed that the group $G=Sz(q)$ is generated by elements $x_1$, $y_1$ and $z_1$ of orders $2, 4$ and $5$ with $x_1y_1z_1=1$, so this gives our first triple $(x_1, y_1, z_1)$. He also showed that $G$ has self-centralising cyclic subgroups of odd orders $q-1$ and $q\pm r+1$, where $r=2^{m+1}$ and $e=2m+1$, and that every element of odd order lies in such a subgroup. Now $q-1$ is coprime to $5$, and either $q+r+1$ or $q-r+1$ is coprime to $5$ as $m\equiv 0$ or $3$ mod~$(4)$ or $m\equiv 1$ or $2$ mod~$(4)$ respectively. We will use this to find a second triple $(x_2, y_2, z_2)$ with elements of orders coprime to those in the first triple.

Taking $G=Sz(q)$ in Proposition~6.1, with $\cal X$ a conjugacy class of elements of order $q-1$, and ${\cal Y}={\cal Z}$ a conjugacy class of elements of order $n=q\pm r+1$, whichever is coprime to $5$, we see from Suzuki's character table of $G$ in~\cite{Suz} that $N({\cal X},{\cal Y},{\cal Z})>0$: every irreducible character $\chi$ takes the value $0$ on either $\cal X$ or $\cal Y$, with the exception of the principal character, taking the value $1$ everywhere, and the character of degree $q^2$, which take the values $1$ and $-1$ on $\cal X$ and $\cal Y$.  Thus $G$ contains a triple $(x_2, y_2, z_2)$ of elements of orders $q-1$, $n$ and $n$ with $x_2y_2z_2=1$.

Suzuki showed that each maximal subgroup of $G$ has order $q^2(q-1)$, $2(q-1)$ or $4(q\pm 2r+1)$, or is isomorphic to $Sz(q')$ where $q'=2^f$ with $e/f$ prime. Simple divisibility arguments show that $x_2$, $y_2$ and $z_2$ cannot be contained in a subgroup of order $q^2(q-1)$, $2(q-1)$ or $4(q\pm 2r+1)$, and by applying Suzuki's classification of the elements of odd order to $Sz(q')$ we see that they cannot be contained in such a subgroup either, so they generate $G$. Since the orders of the elements in this triple are coprime to those in $(x_1, y_1, z_1)$, it follows that $G$ admits a Beauville structure. \hfill$\square$

\medskip

The Ree groups $R(q)={}^2G_2(q)$, introduced by Ree in~\cite{Ree}, are simple groups of order $q^3(q^3+1)(q-1)$, where $q=3^e$ for some odd $e\geq 3$.

\begin{thm}
The Ree group $R(3^e)$ admits a Beauville structure for each odd $e\geq 3$.
\end{thm}

\noindent{\sl Proof.} The argument is similar to that used for the Suzuki groups. In this case we use triples $(x_1, y_1, z_1)$ of orders $2$, $3$ and $7$, discussed by Sah in~\cite{Sah} and by Jones in~\cite{Jon}. We choose $x_2$, $y_2$ and $z_2$ of orders $(q-1)/2$, $n$ and $n$, where $n=q\pm r+1$ with $3r^2=q$, whichever value of $n$ is coprime to $7$. To show the existence of such triples with $x_2y_2z_2=1$ we use the character values given by Ward in~\cite{War}: the only non-principal irreducible character not vanishing at $x_2$ or $y_2$ is that of degree $q^3$, taking the values $1$ and $-1$ respectively. The maximal subgroups of $R(q)$ are given by Levchuk and Nuzhin~\cite{LN} and by Kleidman~\cite{Kle}: they have orders $q^3(q-1)$ or $6(q\pm r+1)$ or $6(q+1)$, or are isomorphic to $R(q')$ where $q'=3^f$ with $e/f$ prime, or to $C_2\times L_2(q)$. It is straightforward to show that $(x_2, y_2, z_2)$ lies in none of these, so this triple generates $R(q)$. Since the orders of $x_1, y_1$ and $z_1$ are coprime to those of $x_2, y_2$ and $z_2$, this shows that $R(q)$ admits a Beauville structure. \hfill$\square$

\medskip

The Beauville structures found here for the Suzuki and Ree groups are not strongly real: there are no automorphisms inverting the elements $y_1$ of orders $4$ and $3$ we have used.

\end{document}